\documentclass{amsart}[12pt]
\usepackage{amssymb}
\usepackage{amsthm}
\usepackage{graphics}
\usepackage{fancyhdr}
\newtheorem*{reftheorem}{Theorem}
\newtheorem{introtheorem}{Theorem}
\newtheorem*{theorem}{Theorem}

\newtheorem{lemma}{Lemma}[section]
\newtheorem{cor}{Corollary}[section]

\newtheorem*{refprop}{Proposition}

\newtheorem*{hamsand}{Ham Sandwich Theorem}

\newtheorem*{boxest}{Box Estimate}

\newtheorem*{bislemma}{Bisection Lemma}
\newtheorem*{vanlemma}{Vanishing Lemma}
\newtheorem*{vislemma}{Visibility Lemma}

\newtheorem*{weakvislemma}{Weak Visibility Lemma}

\numberwithin{equation}{section}
\title{The endpoint case of the Bennett-Carbery-Tao multilinear
Kakeya conjecture}
\author{Larry Guth}
\address{Department of Mathematics, University of Toronto, 40 St.
George St., Toronto ON, Canada}
\email{lguth@math.toronto.edu}

\begin{document}
\begin{abstract} We prove the endpoint case of the multilinear
Kakeya conjecture of Bennett, Carbery, and Tao.  The proof uses the polynomial
method introduced by Dvir.

\end{abstract}

\maketitle

In \cite{BCT}, Bennett, Carbery, and Tao formulated a multilinear
Kakeya conjecture, and they proved the conjecture except for the
endpoint case.  In this paper, we slightly sharpen their result by
proving the endpoint case of the conjecture.

Our method of proof is very different from the proof of Bennett, Carbery,
and Tao.  The original proof was based on monotonicity estimates for heat flows.  
In 2007, Dvir \cite{Dv} made a breakthrough on the Kakeya problem, proving the 
Kakeya conjecture over finite
fields.  His proof used polynomials in a crucial way.  It was not clear whether
Dvir's approach could be adapted to prove estimates in Euclidean space.  Our proof
of the multilinear Kakeya conjecture is based on Dvir's polynomial method.  In my opinion,
the method of proof is as interesting as the result.

The multilinear Kakeya conjecture concerns the overlap properties
of cylindrical tubes in $\mathbb{R}^n$.   Roughly, the (multilinear)
 Kakeya conjecture says that cylinders pointing in different directions cannot
 overlap too much.
 
Before coming to the Bennett-Carbery-Tao multilinear estimate, I want to state a weaker result,
because it's easier to understand and easier to prove.  To be clear about the notation,
a cylinder of radius $R$ around a line $L \subset \mathbb{R}^n$ is the
set of all points $x \in \mathbb{R}^n$ within a distance $R$ of the line
$L$.  We call the line $L$ the core of the cylinder. 

\begin{introtheorem} Suppose we have a finite collection of cylinders $T_{j, a} \subset \mathbb{R}^n$,
where $1 \le j \le n$, and $1 \le a \le A$ for some integer $A$.  Each cylinder has radius 1.
Moreover, each cylinder $T_{j,a}$ runs nearly parallel to the $x_j$-axis.  More precisely, we assume that the angle between the core of $T_{j,a}$ and the $x_j$-axis is at most $(100 n)^{-1}$.

We let $I$ be the set of points that belong to at least one cylinder in each direction.  In symbols,

$$I := \cap_{j=1}^n \left[ \cup_{a=1}^A T_{j,a} \right] . $$

Then $Vol(I) \le C(n) A^{\frac{n}{n-1}}$.

\end{introtheorem}

As Bennett, Carbery, and Tao point out in \cite{BCT}, this
estimate can be viewed as a generalization of the Loomis-Whitney
inequality.

\begin{reftheorem}(special case of Loomis and Whitney, 1949, \cite{LW}) Let $U$
be an open set in $\mathbb{R}^n$.  Let $\pi_j$ denote the
projection from $\mathbb{R}^n$ onto the hyperplane perpendicular
to the $x_j$-axis.  Suppose that for each $j$, $\pi_j(U)$ has (n-1)-dimensional
volume at most $B$.

Then $Vol(U) \le B^{\frac{n}{n-1}}$.

\end{reftheorem}

Suppose that each tube $T_{j,a}$ runs exactly parallel to the
$x_j$-axis.  It follows that $\pi_j(I)$ is contained in $A$
unit balls and has volume at most $\omega_{n-1} A$.  Applying
the Loomis-Whitney inequality, we see that the volume of $I$ is
bounded by $\lesssim A^\frac{n}{n-1}$.   Theorem 1 says that - 
up to a constant factor - this volume 
estimate continues to hold if we allow the tubes to tilt slightly.  

The Loomis-Whitney inequality is sharp whenever the open set $U$
is a cube.  Similarly, Theorem 1 is essentially sharp
whenever the tubes are arranged in a cubical lattice.

The proof of Theorem 1 uses the polynomial method of Dvir.  The main
new idea in the paper is a new approach for adapting Dvir's method
to $\mathbb{R}^n$.  The new approach uses algebraic topology.  In
particular, we will use a polynomial generalization of the ham sandwich
theorem, proven by Stone and Tukey \cite{ST} in the early 40's.

Now we turn to the multilinear version of the Kakeya maximal conjecture,
formulated by Bennett, Carbery, and Tao.

\begin{introtheorem} (Multilinear Kakeya estimate)

For each $1 \le j \le n$, let $T_{j,a}$ be a collection of unit cylinders,
where $a$ runs from $1$ to $A(j)$.  We let $v_{j,a}$ be a unit vector parallel to
the core of $T_{j,a}$.  We assume that the cylinders from different classes are quantitatively
transverse in the sense that any determinant of a matrix $(v_{1, a_1}, v_{2, a_2}, ... , v_{n,
a_n})$ has norm at least $\theta > 0$. 

Under these hypotheses, the following inequality holds.

$$ \int \left[ \prod_{j=1}^n \left( \sum_{a=1}^{A(j)} \chi_{T_{j,a}} \right)
\right]^{\frac{1}{n-1}} \le C(n) \theta^{\frac{-1}{n-1}} \prod_{j=1}^n
A(j)^{\frac{1}{n-1}}.$$

\end{introtheorem}

Theorem 2 generalizes Theorem 1.  If each vector $v_{j,a}$ lies within a 
small angle of the $x_j$-axis, then the determinant condition is easy to check, and so Theorem 2 applies.  Recall
that $I$ is the set of points lying in at least one cylinder with each value of $j$.  At every point $x \in I$,
the integrand in Theorem 2 is at least 1.  Hence Theorem 2 gives an upper bound for
the volume of $I$, recovering Theorem 1.

Theorem 2 improves on Theorem 1 in the following ways.  First, we allow a more general
condition on the angles of the tubes.  Second, we allow the different classes to have different numbers
of tubes: $A(j)$ depends on $j$.  Third, and most importantly, we get an integral bound where the
integrand is very large at ``high multiplicity" points - points which lie in many tubes from each
direction.

The paper \cite{BCT} has a very nice introductory discussion of
the multilinear Kakeya estimate.  Some of the topics it describes
are the original Kakeya conjecture, and linear and multilinear
restriction estimates.  Using their multilinear Kakeya estimates,
Bennett, Carbery, and Tao are able to prove nearly optimal
multilinear restriction estimates.  We refer to that paper for
more context.

The proof of Theorem 2 is harder than the proof of Theorem 1.
It uses more sophisticated tools from algebraic topology:
cohomology classes, cup products, and the Lusternik-Schnirelmann vanishing
theorem.  Theorem 2 is more important than Theorem 1 because Bennett, Carbery,
and Tao use Theorem 2 to prove $L^p$ estimates for multilinear restriction
operators.  On the other hand, Theorem 1 contains the main ideas of this paper, 
and its proof is only three pages long.

This paper uses algebraic topology.  I want it to be
understandable to mathematicians who work in analysis and combinatorics,
so I will try to introduce the algebraic topology in a friendly way.  In particular,
there is a short section introducing Lusternik-Schnirelmann theory, and an
appendix giving the proof of the Lusternik-Schnirelmann vanishing lemma.

As a corollary of our method, we give a `planiness' estimate for
unions of tubes in $\mathbb{R}^n$.  An estimate of this kind can
also be proven using the methods of \cite{BCT}, but the estimate
below is slightly sharper.  The phenomenon of `planiness' was
discovered by Katz, Laba, and Tao in \cite{KLT}, and the estimate
below is similar to some estimates from that paper.

\begin{boxest} There is a constant $C(n) > 0$ so that the
following holds.  Suppose $X \subset \mathbb{R}^n$ is a union of
cylinders with radius 1 and length $L >> 1$.  For each $x \in X$
we can choose a rectangular box $B(x)$ with the following
properties.

1. The box $B(x)$ is centered at $x$.  It may be oriented in any
direction.  It has volume at most $C(n) Vol(X)$.

2. For every cylinder $T \subset X$ of radius 1 and length $L$,
if we pick a random point $x \in T$, then with probability at
least $9/10$, the tube $T$ lies in the box $B(x)$.

\end{boxest}

\vskip5pt

{\bf Acknowledgements.} I would like to thank Nets Katz for
showing me the multilinear Kakeya estimates in
\cite{BCT}.  I showed him the proof of the box estimate, and he
explained to me how that estimate is related to multilinear
Kakeya estimates and the work of Bennett, Carbery, and Tao.  I would also like to
thank Kannan Soundararajan for interesting conversations about
combinatorial number theory.  In particular, he pointed out to me
Dvir's paper \cite{Dv}.

\section{The polynomial ham sandwich theorem}

The main tool in our proof is a generalization of the ham
sandwich theorem to algebraic hypersurfaces.  I learned about this
result from Gromov's paper \cite{Gr}.  However, I recently learned
that it was proven by Stone and Tukey \cite{ST} in 1941.
In this section, we explain and prove this generalization of the ham sandwich theorem, 
following Stone and Tukey.

First we recall the
original ham sandwich theorem.

\begin{hamsand} Let $U_1$, ..., $U_n$ be finite volume open
sets in $\mathbb{R}^n$.  Then there is a hyperplane $H$ that
bisects each set $U_i$.
\end{hamsand}

The 3-dimensional case of the ham sandwich theorem was first
proven in the 30's by Stefan Banach, using the Borsuk-Ulam
theorem.  Stone and Tukey extended the method to the n-dimensional
case.  (There is a nice historical discussion Wikipedia.)  Stone and Tukey 
noticed that the same method can be used to prove many other bisection
results.  For example, they proved the following
proposition.

\newtheorem*{polyham}{Polynomial Ham Sandwich Theorem}

\begin{polyham} (Stone, Tukey \cite{ST})
Let $N = {{n+d} \choose d} - 1$.  Let $U_1$, ..., $U_N$ be finite
volume open sets in $\mathbb{R}^n$.  Then there is a degree d algebraic
hypersurface $Z$ which bisects each set $U_i$.
\end{polyham}

We will prove the polynomial ham sandwich
theorem using the Borsuk-Ulam theorem.  We recall the Borsuk-Ulam theorem.

\newtheorem*{BU}{Borsuk-Ulam Theorem}

\begin{BU}  Let $F$ be a continuous map from $S^N$ to
$\mathbb{R}^N$ obeying the antipodal condition

$$ F(-x) = - F(x)  \textrm{  for every $x \in S^N$. }  $$

Then the image of $F$ contains 0.

\end{BU}

For a proof of the Borsuk-Ulam theorem, the reader may look at Hatcher's
book on algebraic topology \cite{H}, pages 174-176.  Another reference is
{\it Using the Borsuk-Ulam Theorem} \cite{M} by Matou\v{s}ek.  This book
gives a proof of the theorem, and it also discusses interesting applications
of the Borsuk-Ulam theorem, for example to Kneser's conjecture in combinatorics.
We now turn to the proof of the polynomial ham sandwich theorem.

\begin{proof} Let $V(d)$ denote the vector space of all real polynomials of
degree at most $d$ in $n$ variables.  The dimension of $V(d)$ is
${{n+d} \choose d}$.  Let $S^N$ denote the unit sphere in $V(d)$,
where recall $N = {{n+d} \choose d} -1$.  For each set $U_i$, we
define a function $F_i$ from $S^N$ to $\mathbb{R}$, by setting

$$F_i(P) = Vol ( \{ x \in U_i | P(x) > 0 \} ) - Vol ( \{ x \in U_i | P(x) < 0 \} ) . $$

If we replace $P$ with $-P$, then the two volumes trade places, so $F_i( - P) =
- F_i (P)$.  It's not hard to check that $F_i$ is continuous (see below for the details).
Combining all $F_i$ into a vector-valued function, we get a continuous map 
$F: S^N \rightarrow \mathbb{R}^N$ obeying the antipodal condition.  By the
Borsuk-Ulam theorem, $F(P) = 0$ for some $P \in S^N \subset V(d)$.  By the definition of
$F_i$, we see that for each $i$

$$Vol ( \{ x \in U_i | P(x) > 0 \} ) = Vol ( \{ x \in U_i | P(x) < 0 \} ) . $$

Hence the hypersurface defined by $P(x) = 0$ bisects each set $U_i$. \end{proof}

For the sake of completeness, we include the proof that $F_i$ is a continuous
function.

\newtheorem*{contlemma}{Continuity Lemma}
\begin{contlemma}  If $U$ is an open set of finite measure, then
the measure of the set $\{ x \in U | P(x) > 0 \}$ depends
continuously on $P \in V(d) \backslash 0$.
\end{contlemma}

\proof Suppose that $P$ is a non-zero polynomial in $V(d)$ and
$P_n \in V(d)$ with $P_n \rightarrow P$.  Pick any $\epsilon >
0$.  We can find a subset $E \subset U$ so that $P_n \rightarrow P$ 
uniformly pointwise on $U - E$, and $m(E) < \epsilon$.  

The set $\{ x \in U | P(x) = 0 \}$ has measure zero.  Therefore, we can
choose $\delta$ so that the set $\{x \in U \textrm{ such that }
|P(x)| < \delta \}$ has measure less than $\epsilon$.  

Next we
choose $n$ large enough so that $|P_n(x) - P(x)| < \delta$ on $U - E$.  Then
the measures of $\{ x \in U | P_n(x) > 0 \}$ and $\{ x \in U |
P(x) > 0 \}$ differ by at most $2 \epsilon$.  But $\epsilon$ was
arbitrary. \endproof

To make use of the polynomial ham sandwich theorem, we will use
a standard volume estimate for hypersurfaces that bisect simple sets.

\newtheorem*{area}{Basic Area Estimate}

\begin{area} If a hypersurface $S$ bisects a unit ball or a unit cube, then $S$
has $(n-1)$-dimensional volume at least $c(n)$.
\end{area}

\section{Directed volume}

The second tool in our paper is directed volume, which is a way
of measuring the amount of volume of a hypersurface facing in
different directions.

For a hypersurface $S \subset \mathbb{R}^n$, we define
a directed volume function $V_S$ by the following formula.

$$V_S(v) := \int_S |v \cdot N| dvol_S. \eqno{(1)} $$

In this formula, $N$ denotes the normal vector to $S$, and
$v \in \mathbb{R}^n$ is a fixed vector.  Hence the directed
volume is a non-negative function of $ v \in \mathbb{R}^n$.

For a unit vector $v$, the directed volume $V_S(v)$ can be
given a different, more geometric interpretation.  Let
$\pi_v: \mathbb{R}^n \rightarrow v^\perp$ be the orthogonal
projection onto $v^\perp$.  Then we can also think of $V_S(v)$ as
the volume of $\pi_v(S)$, counted with geometric multiplicity.  
For each $y \in v^\perp$, we consider
the intersection $S \cap \pi_v^{-1}(y)$.  We let $|S \cap \pi_v^{-1}(y)|$
denote the number of points in $S \cap \pi_v^{-1}(y)$.  For a compact
smooth hypersurface $S$ (possibly with boundary), this number of points
is finite for almost every $y$.  If $v$ is a unit vector, then
$V_S(v)$ is given by the following formula.

$$V_S(v) = \int_{v^\perp} |S \cap \pi_v^{-1}(y)| dy. \eqno{(2)}  $$

Equations (1) and (2) will both be useful to us.  Using equation
(2), we can prove a key estimate about the directed volumes of
algebraic hypersurfaces in cylinders.

\begin{lemma} (Cylinder estimate) If $T$ is a cylinder of radius $r$,
$v$ is a unit vector parallel to the core of $T$, and $Z$ is an
algebraic hypersurface of degree $d$, then the directed volume
$V_{Z \cap T}(v)$ is bounded as follows:

$$ V_{Z \cap T}(v) \le \omega_{n-1} r^{n-1} d. $$

\end{lemma}

\proof The projection $\pi_v(T)$ is an (n-1)-dimensional disk of
radius $r$.
The function $|Z \cap T \cap \pi_v^{-1}(y)|$ is supported in this
disk.  But since $Z$ is a degree $d$ algebraic hypersurface, $Z$
intersects a line in at most $d$ points, unless $Z$ contains the
entire line.  Hence $|Z \cap \pi_v^{-1}(y)| \le d$ for almost every
$y$. \endproof

Equation (1) is also useful.  For example, it allows us to see
that a surface of volume 1 must have a fairly large directed
volume in some direction.

\begin{lemma} Suppose that $v_1$, ..., $v_n$ are unit vectors.
Let $e_j$ denote the coordinate unit vectors, and suppose that
$|e_j - v_j| < (100n)^{-1}$.  Let $S$ be any hypersurface in
$\mathbb{R}^n$.  Then $Vol(S) \le 2 \sum_{j=1}^n V_S(v_j)$.
\end{lemma}

\begin{proof} For each $x$ in $S$, let $N(x)$ denote the unit normal vector
to $S$ at $x$.  Because of the angle condition on $v_j$, we
know that $| v_j \cdot N(x)| \ge |e_j \cdot N(x)| - (100 n)^{-1}$.
Hence $\sum_j  |v_j \cdot N(x)| \ge (\sum_j |e_j \cdot N(x)|) - (100)^{-1}
\ge (99/100)$.

Integrating this inequality over $S$, we see that

$$ \sum_j  V_S(v_j) = \int_S \sum_j |v_j \cdot N(x)| \hskip3pt dvol_S (x) \ge \int_S \frac{99}{100} dvol_S(x) = \frac{99}{100} Vol(S). $$

\end{proof}

Estimates for directional volumes appeared in some papers that I
wrote giving quantitative estimates for certain homotopy
invariants of a map in terms of its Lipschitz constant -
\cite{G2} and \cite{G3}.

\section{The proof of Theorem 1}

In this section, we prove Theorem 1.

\begin{theorem} Suppose we have a finite collection of cylinders $T_{j, a} \subset \mathbb{R}^n$,
where $1 \le j \le n$, and $1 \le a \le A$ for some integer $A$.  Each cylinder has radius 1.
Moreover, each cylinder $T_{j,a}$ runs nearly parallel to the $x_j$-axis.  More precisely, we assume that the angle between the core of $T_{j,a}$ and the $x_j$-axis is at most $(100 n)^{-1}$.

We let $I$ be the set of points that belong to at least one cylinder in each direction.  In symbols,

$$I := \cap_{j=1}^n \left[ \cup_{a=1}^A T_{j,a} \right] . $$

Then $Vol(I) \le C(n) A^{\frac{n}{n-1}}$.

\end{theorem}

\proof Look at the standard unit lattice in $\mathbb{R}^n$.  Let
$Q_1, ..., Q_V$ be the set of n-cubes in the lattice which
intersect $I$.  Here $V$ is the number of cubes that intersect
$I$.  It suffices to prove the estimate $V \lesssim
A^{\frac{n}{n-1}}$.

According to the polynomial ham sandwich theorem, we may find a degree $d$ algebraic
hypersurface $Z$ which bisects $Q_k$ for every $k$, with degree
$d \lesssim V^{1/n}$.  Because of the bisection property, the
volume of $Q_k \cap Z$ is $\gtrsim 1$ for each $Q_k$.

For each $Q_k$, we pick a tube in each direction that goes
through $Q_k$.  So we have labels $a_1(k), ..., a_n(k)$ so that
$T_{j, a_j(k)}$ intersects $Q_k$.  By assumption, the vector
$v_{j, a_j(k)}$ is within $(100n)^{-1}$ of the coordinate vector
$e_j$.  Applying Lemma 2.2, we get the following estimate.

$$\sum_{j=1}^n V_{Z \cap Q_k} (v_{j, a_j(k)}) \gtrsim Vol(Z \cap Q_k) \gtrsim 1.$$

So for each $k$, we can choose a tube $T_{j(k), a(k)}$ which
meets $Q_k$ and so that $V_{Z \cap Q_k} (v_{j(k), a(k)}) \gtrsim
1$.

We have just associated a tube to each cube.  There are in total
only $n A$ tubes.  By the pigeonhole principle, there is a tube
associated to $\gtrsim V/A$ different cubes. Let this tube be
$T_{j,a}$.  Then we have $\gtrsim V/A$ different cubes $Q_k$
which intersect $T_{j,a}$ and with $V_{Z \cap Q_k}(v_{j,a})
\gtrsim 1$.

Let $\tilde T_{j,a}$ denote the $\sqrt n$ neighborhood of
$T_{j,a}$.  The set $\tilde T_{j,a}$ is itself a cylinder of
radius $1 + \sqrt n$, with core parallel to $v_{j,a}$,
and it contains all the cubes $Q_k$ which overlap $T_{j,a}$.

Therefore, the directed volume $V_{Z \cap \tilde T_{j,a}}
(v_{j,a}) \gtrsim V/A$.

On the other hand, by the cylinder estimate in Lemma 2.1, the same directed
volume is $\lesssim V^{1/n}$.

Hence $V/A \lesssim V^{1/n}$.  Rearranging, we get $V \lesssim
A^{\frac{n}{n-1}}$. \endproof

\section{The Lusternik-Schnirelmann vanishing lemma}

To prove Theorem 2, we use some more sophisticated algebraic topology:
the Lusternik-Schnirelmann vanishing lemma.  In this section, I will
introduce it and try to explain what it's good for.  The basic message is that
the vanishing lemma is similar to the ham sandwich theorem, but it's more
flexible.

The vanishing lemma is about cup-products of cohomology classes.

\begin{vanlemma} Let $X$ be a CW complex (for example a manifold).  Let $a_1, a_2$ be cohomology classes in
$H^*(X, R)$, where $R$ may be any ring of coefficients, such as $\mathbb{R}$, $\mathbb{Z}$, or
$\mathbb{Z}_2$.  Suppose that $a_1 $ vanishes on some open set $S_1 \subset X$ and that $a_2$ vanishes on some open set $S_2 \subset X$.  Then the cup product $a_1 \cup
a_2$ vanishes on the union $S_1 \cup S_2$.
\end{vanlemma}

The vanishing lemma is one of the fundamental topological facts about cup products.  I believe that it
was first proven by Lusternik and Schnirelmann in the 1930's, as part of their project for proving the existence of closed geodesics.  The proofs I have seen in the literature are a little more abstract than I would like, so I wrote an appendix giving the proof.

Here is the basic intuition behind the vanishing lemma.  Suppose that $f_1$ and $f_2$ are functions on $X$.  If $f_1$ vanishes on $S_1$ and $f_2$ vanishes on $S_2$, then clearly the product $f_1 f_2$ vanishes
on the union $S_1 \cup S_2$.  The vanishing lemma holds because cohomology classes are not so
different from functions.  A cohomology class can be represented by either a differential form or a singular cocycle, and these objects have enough in common with functions to make the vanishing lemma hold.  For details, see the appendix.

To apply the vanishing lemma, we need to know something about the cup products of cohomology classes.  For this paper, the key example is the cohomology ring of real projective space.

\newtheorem*{cohr}{Cohomology ring of $\mathbb{RP}^N$}

\begin{cohr} The cohomology group $H^i (\mathbb{RP}^N, \mathbb{Z}_2)$ is isomorphic to $\mathbb{Z}_2$ for $0 \le i \le N$, and is equal to 0 otherwise.  Let $a$ denote the non-zero element in $H^1(\mathbb{RP}^N, \mathbb{Z}_2)$.  Then for $1 \le i \le N$, $a^i$ is the non-zero element of
$H^i(\mathbb{RP}^N, \mathbb{Z}_2)$.
\end{cohr}

This theorem may be found in Hatcher's topology book \cite{H} on page 212.

Using the vanishing lemma, we can give a different proof of the polynomial
ham sandwich theorem.

As before, we let $V(d)$ denote the vector space of all real polynomials of
degree at most $d$ in $n$ variables.  The dimension of $V(d)$ is
${{n+d} \choose d}$.  For each non-zero polynomial $P$ in $V(d)$,
there is an associated variety, the zero-set of $P$.  If we
replace $P$ by some multiple $\lambda P$, the zero-set remains
unchanged, and so the real algebraic hypersurfaces of degree at
most $d$ are parametrized by the projectivization of $V(d)$, which
is a real projective space $\mathbb{RP}^{N}$, where $N = {{n+d}
\choose d} - 1$. 

We're interested in hypersurfaces that bisect open sets.  Given a finite
volume open set $U \subset \mathbb{R}^n$, we let $Bi(U) \subset \mathbb{RP}^N$
consist of the algebraic hypersurfaces that bisect the set $U$.  If $Z$ is
a real algebraic hypersurface given by the equation $P = 0$, then we say
that $Z$ bisects $U$ if

$$Vol \{ x \in U | P(x) > 0 \} = Vol \{x \in U | P(x) < 0 \}.$$

By the continuity lemma from Section 1, these volumes change continuously
with $P$, and so $Bi(U)$ is a closed subset of $\mathbb{RP}^N$.  The key
topological result about $Bi(U)$ is the following lemma.

\begin{bislemma} Let $a$ denote the non-trivial cohomology class
in $H^1(\mathbb{RP}^N, \mathbb{Z}_2)$.  Let $U$ be a finite-volume subset of
$\mathbb{R}^n$.  Then the
cohomology class $a$ vanishes on the complement $\mathbb{RP}^N -
Bi(U)$.
\end{bislemma}

\proof Suppose that $a$ does not vanish on $\mathbb{RP}^N -
Bi(U)$.  Then the class $a$ is detected by a 1-cycle $c$ in
$\mathbb{RP}^N - Bi(U)$.  Without loss of generality, we may
assume that $c$ has only one component, and so $c$ is
topologically a circle. 
Pick a point in $c$ and look at the corresponding hypersurface $Z$. 
We can assume that $Z$ does not bisect $U$, so the complement of
$Z$ has a big half and a little half.  Now we pick a polynomial
$P_Z$ representing $Z$, and we choose it so that the big half is
where the polynomial $P_Z$ is positive.  We can lift our
1-parameter family of hypersurfaces to a 1-parameter family of
polynomials that goes from $P_Z$ to $-P_Z$.  The part of $U$
where $P_Z$ is positive has more than half measure.  The part of
$U$ where $-P_Z$ is positive has less than half measure.  According
to the Continuity Lemma from Section 1, the measure changes continuously
as the polynomial changes.  By
continuity, there is a polynomial in the family that bisects $U$.
\endproof

Combining the bisection lemma and the vanishing lemma, we 
can say something about hypersurfaces that bisect multiple
sets.  Suppose that $U_1, ..., U_r \subset \mathbb{R}^n$ are finite volume open sets,
where $r$ is any positive integer.  Let $Bi(U_1, ..., U_r)
\subset \mathbb{RP}^N$ denote the set of algebraic hypersurfaces
which bisect all the open sets $U_1, ..., U_r$.  The set $Bi(U_1, ..., U_r)$
is just the intersection of $Bi(U_i)$ ($1 \le i \le r$).  In particular, $Bi(U_1, ..., U_r)
\subset \mathbb{RP}^N$ is a closed set.

\newtheorem*{multbis}{Multiple bisection lemma}
\begin{multbis}  Let $Bi(U_1, ..., U_r)$ be as above. Then the cohomology
class $a^r$ vanishes on $\mathbb{RP}^N - Bi(U_1, ..., U_r)$.
\end{multbis}

\proof By the bisection lemma, the cohomology class $a$ vanishes
on $\mathbb{RP}^N - Bi(U_i)$ for each $i$.  Each of these sets is
open.  Therefore, the vanishing lemma tells us that $a^r$ vanishes
on their union.  But the union $\cup_{i=1}^r [ \mathbb{RP}^N - Bi(U_i) ]$
is exactly $\mathbb{RP}^N - Bi(U_1, ..., U_r)$. \endproof

Combining the multiple bisection lemma and the cohomology ring of
$\mathbb{RP}^N$, we can reprove the polynomial ham sandwich
theorem.  This proof was given by Gromov in \cite{Gr}.

\begin{polyham} Let $U_1, ..., U_N$ be any finite volume subsets of
$\mathbb{R}^n$, where $N = {{n+d} \choose d} - 1$.  
Then there is a real algebraic
hypersurface of degree at most $d$ that bisects each set $U_i$.
\end{polyham}

\proof Recall that $a$ is the non-zero cohomology class in
$H^1(\mathbb{RP}^N, \mathbb{Z}_2)$.  By the multiple bisection lemma, $a^N$
vanishes on $\mathbb{RP}^N - Bi(U_1, ..., U_N)$.  But in the cohomology
ring of $\mathbb{RP}^N$, $a^N$ does not vanish on $\mathbb{RP}^N$.
Hence $Bi(U_1, ..., U_N)$ must be non-empty.  In other words,
there is a degree $d$ hypersurface $Z$ that bisects each open
set $U_i$.  \endproof

The vanishing lemma has other applications besides the ham sandwich theorem.  One classical application is to give covering estimates.

\newtheorem*{covest}{Covering Estimate}

\begin{covest} (Lusternik-Schnirelmann) Suppose that $\mathbb{RP}^N$ is covered by some contractible open sets $S_1$, ..., $S_r$.  Then $r \ge N+1$.
\end{covest}

\begin{proof} Since each $S_i$ is contractible, the cohomology class $a$ vanishes on each $S_i$.  Applying the vanishing lemma once, we see that $a^2$ vanishes on $S_1 \cup S_2$.  Proceeding
inductively, we see that $a^r$ vanishes on the union of all $S_r$, which is $\mathbb{RP}^N$.  But
in the cohomology ring of $\mathbb{RP}^N$, $a^i$ is non-zero for all $i \le N$.  Hence $r \ge N+1$.
\end{proof}

Our proof of the multilinear Kakeya estimate combines the polynomial ham sandwich theorem with some covering estimates similar to the one above.  The Lusternik-Schnirelmann vanishing lemma allows
us to combine these two techniques, making it a little more flexible than the Borsuk-Ulam theorem.

\section{The visibility lemma}

In the proof of Theorem 1, we found an algebraic hypersurface whose
intersection with many unit cubes has volume $\gtrsim 1$.  The total
volume is not as important as the directional volumes $V_{Z \cap Q}(v)$
in various directions.  (Here $Z$ is the hypersurface, $Q$ is a cube,
and $v$ is a direction.)  In this section, we build a hypersurface
which has large directional volumes in many directions.

We now define the visibility of a hypersurface, which measures whether
the surface has a large directional volume in many directions.
Roughly, a hypersurface has large visibility
if either $V_S(v)$ is large for every unit vector $v$, or else
$V_S(v)$ is extremely large for some vectors $v$.

We define the ``visibility'' of a surface $S$ to be

$$Vis[S] := Vol \left( \{ v \textrm{ such that } |v| \le 1 \textrm{ and }
V_S(v) \le 1 \} \right) ^{-1}.$$

This definition is a little long, so we make some comments about it.  The
reader may wonder, why not just look at the average directional volume in all unit
directions $v$: $  \int_{S^{n-1}} V_S(v) dvol(v) / Volume(S^{n-1})$?  For our
arguments, it's crucial to know whether $V_S(v)$ is small in some directions even
if the set of such directions has a small measure.  The average directional volume
above won't detect small values of $V_S(v)$, but the definition of visibility is
quite sensitive to small values of $V_S(v)$.

We now compute the visibility in two examples. 
First, suppose that $S$ is a unit (n-1)-disk in the hyperplane $x_n = 0$.  
Then the function $V_S(v) =
\omega_{n-1} |v_n|$, where $v_n$ is the $n^{th}$ component of
$v$.  Therefore, the set $\{ v | V_s(v) \le 1 \}$ is an infinite
slab of the form $|v_n| \le C$.  The set of $v$ with $V_S(v) \le
1$ and $|v| \le 1$ is roughly the unit ball, and so the
visibility of $S$ is on the order of 1.  We had to include the
condition $|v| \le 1$ in the definition, or else the visibility
of a disk would be zero.  Including $|v| \le 1$ in the definition
has the unpleasant effect that the visibility of the empty set is
also around 1.  In practice, we will speak of the visibility of
$S$ for surfaces $S$ contained in a unit cube and with volume at
least 1, and in this range, the visibility behaves reasonably. 

A second important example is a union of unit disks with $N_j$ disks
perpendicular to the $x_j$-axis.  In this case, $V_S(v)$ is roughly
$\sum_j N_j |v_j|$.  Hence the region where $V_S(v) \le 1$ is
roughly $\{ v \in \mathbb{R}^n | |v_j| \le N_j^{-1} \}$.  The volume of this
region is roughly $N_1^{-1} ... N_n^{-1}$, and so the visibility of this surface
is roughly $N_1 ... N_n$.  This is the best example to keep in mind to
understand what visibility means.

Our next goal is to find algebraic hypersurfaces which have large
visibility in many cubes.  We recall from
Section 4 that the space of degree $d$ hypersurfaces in
$\mathbb{R}^n$ is parametrized by $\mathbb{RP}^N$ for $N = {{n+d}
\choose d} - 1$.  We will slightly abuse notation by identifying
a degree $d$ hypersurface $Z$ and the corresponding point in
$\mathbb{RP}^N$ - we will speak of $Z \in \mathbb{RP}^N$.

If we fix some cube $Q \subset \mathbb{R}^n$, we want to study
$Vis [ Z \cap Q]$ as a function of $Z \in \mathbb{RP}^N$.  Unfortunately,
this function is not continuous in $Z$.  Even the (n-1)-dimensional volume of
$Z \cap Q$ is not continuous in $Z$.  Because we make
topological arguments using the Lusternik-Schnirelmann vanishing
lemma, this discontinuity leads to some technical problems.  To
deal with these, we define mollified continuous versions of the
directed volume and the visibility.  We mollify these functions by averaging over small balls in $\mathbb{RP}^N$.  We use the
standard metric on $\mathbb{RP}^N$, and let $B(Z, \epsilon)$
denote the ball around $Z \in \mathbb{RP}^N$ of radius $\epsilon$.

For any open
set $U$, we define a mollifed version of $V_{Z \cap U}(v)$ as
follows.

$$\bar V_{Z \cap U} (v) := |B(Z, \epsilon)|^{-1}
\int_{B(Z, \epsilon)} V_{Z' \cap U}(v) dZ'.$$

We define a mollified visibility function using the mollified directional volumes.

$$\overline{Vis}[Z \cap U] := Vol \left( \{ v \textrm{ such that } |v| \le 1
\textrm{ and } \bar V_{Z \cap U} (v) \le 1 \} \right)^{-1}.$$

We will choose $\epsilon$ extremely small compared to all other
constants in the paper.  In practice, the mollified directional
volumes and visibilities maintain all the useful properties of
the unmollified versions, and they are also continuous. 
Therefore, on an early reading of the paper, it makes sense to
ignore the mollification and just pretend that the visibility is
continuous in $Z$.

In the following lemma, we collect the properties of the
mollified volumes and visibilities which we will use.

\begin{lemma} Let $U$ be a bounded open set in $\mathbb{R}^n$. 
The mollified directed volume $\bar V_{Z \cap
U}(v)$ and the mollified visibility $\overline{Vis}[Z \cap U]$ obey
the following properties.

(i) Scaling: For any constant $\lambda$, $\bar V_{Z \cap U}
(\lambda v) = |\lambda| \bar V_{Z \cap U}(v)$.

(ii) Convexity: The function $\bar V_{Z \cap U}(v)$ is convex
in $v$.

(iii) Disjoint unions: If $U$ is a disjoint union of $U_1$ and
$U_2$, then $\bar V_{Z \cap U}(v) = \bar V_{Z \cap U_1}(v) +
\bar V_{Z \cap U_2}(v)$.

(iv) Cylinder estimate: If $T$ is a cylinder of radius $r$
with core vector $v$, and if $Z$ is a degree $d$ hypersurface,
then $\bar V_{Z \cap T}(v) \le \omega_{n-1} r^{n-1} d$.

(v) Bisection: If $Z$ bisects a unit ball $B$, and if $\epsilon$
is small enough, then $\bar V_{Z \cap B}(v) \gtrsim 1$ for some
unit vector $v$.

(vi) Continuity: The functions $\bar V_{Z \cap U}(v)$ and 
$\overline{Vis}[Z \cap U]$ depend continuously on $Z \in \mathbb{RP}^N$.

\end{lemma}

\proof (i) This follows by plugging in the formulas.

(ii) For each vector $N$, the function $|N \cdot v|$ is a convex
function of $v$.  Since a positive combination of convex
functions is convex, $V_{Z \cap U}(v)$ is convex in $v$.  Since an
average of convex functions is convex, $\bar V_{Z \cap U}$ is
also convex.

(iii) This also follows by plugging in the formulas.

(iv) Lemma 2.1 tells us that $V_{Z' \cap T}(v) \le \omega_{n-1}
r^{n-1} d$ for every degree d hypersurface $Z'$.  Taking an
appropriate average, we see that $\bar V_{Z \cap T}(v) \le
\omega_{n-1} r^{n-1} d$.

(v) Suppose that $Z$ bisects $B$.  By the Continuity Lemma in Section 1,
we can choose $\epsilon$ small enough so that each $Z'$ in $B(Z,
\epsilon)$ nearly bisects $B$.  Hence the volume of $Z' \cap B$
is $\gtrsim 1$.  We let $e_1, ..., e_n$ be the standard
orthonormal basis of $\mathbb{R}^n$.  By Lemma 2.2, $\sum_{j=1}^n
V_{Z' \cap B}(e_i) \ge (1/2) Vol(Z' \cap B) \gtrsim 1$.  Taking an
average over $Z'$ in $B(Z, \epsilon)$, we see that $\sum_{j=1}^n
\bar V_{Z \cap B}(e_i) \gtrsim 1$.

(vi) The function $V_{Z \cap U}(v)$ is a bounded measurable
function on $\mathbb{RP}^N$.  Hence its averages over
$\epsilon$-balls form a continuous function.  So $\bar V_{Z
\cap U}(v)$ depends continuously on $Z$.  

Next we address continuity in $v$.  The function $V_{Z \cap U}(v)$ is
Lipschitz in $v$ with a constant $C(d, U, n)$ independent
of $Z$.  To see this, we expand $|V_{Z\cap U}(v_1) - V_{Z \cap U}(v_2)|$
as an integral : $| \int_{Z \cap U} |N \cdot (v_1 -v_2) | dvol | \le |v_1 - v_2| Vol(Z \cap U)$.
Now $U$ is a bounded domain, so it fits in a ball of some radius $R(U)$, and
standard algebraic geometry shows that $Vol(Z \cap U) \le C_n R^{n-1} d$.  (For
more details on this Crofton estimate, see \cite{G} page 58.)

Hence $\bar V_{Z \cap U}(v)$ is also Lipschitz in $v$ with a constant $C(d, U,n)$.  Therefore,
$\bar V_{Z \cap U}(v)$ is jointly continuous as a function of $(Z,v) \in \mathbb{RP}^N \times \mathbb{R}^n$.  Hence $\overline{Vis}[Z \cap U]$ is continuous in
$Z$. \endproof

The next lemma allows us to find algebraic hypersurfaces
with large visibility.  The lemma is analogous to the bisection lemma,
but instead of producing surfaces that bisect a ball, it produces surfaces
with large visibility in a ball.

\begin{vislemma} There is an integer constant $C_n > 1$ so that the
following holds.  Fix any degree $d$ and any unit ball $B(p,1)
\subset \mathbb{R}^n$.  Consider the space of degree $d$
algebraic hypersurfaces in $\mathbb{R}^n$, parametrized by
$\mathbb{RP}^{N}$.  Let $L_M$ denote the subset of
algebraic surfaces $Z$ with $\overline{Vis}[Z \cap B(p,1)] \le M$,
where $M \ge 1$ is an integer.  Let $a$ denote the non-zero
cohomology class in $H^1(\mathbb{RP}^{N}, \mathbb{Z}_2)$. 
Then the cohomology class $a^{C_n M}$ vanishes on a
neighborhood of $L_M$.
\end{vislemma}

\proof Let $E$ be an ellipsoid contained in the unit ball in
$\mathbb{R}^n$, with the volume of $E$ at least $M^{-1}$.  Let
$L(E)$ denote the set of degree d hypersurfaces $Z$ so that

$$\bar V_{Z \cap B(p,1)} (v) \le 1, \textrm{  for all } v \in E.$$

Notice that if $\overline{Vis}[Z \cap B(p,1)] \le M$, then $Z$
is in $L(E)$ for some ellipsoid $E$ of volume $\gtrsim
M^{-1}$.  We will first deal with the different ellipsoids
individually and then see how to deal with all of them
simultaneously.

\begin{weakvislemma} If $E$ is an ellipsoid contained in the unit
ball with volume at least $M^{-1}$, then the cohomology class
$a^{C_n M}$ vanishes on a neighborhood of $L(E)$.
\end{weakvislemma}

\proof Let $A(n)$ be a large number we will choose later. 

We let $E'$ be a rescaling of $E$ by a factor $A(n)^{-1}$ (so
that $E'$ is smaller than $E$).  We let $U_1$, ..., $U_k$ denote
disjoint parallel copies of $E'$ contained in $B(p,1)$.  We take 
a maximal family of parallel copies of $E'$ in $B(p,1)$ - meaning
that there is not room to add an additional parallel copy of $E'$.
From the maximality, we see that $Vol(E') k \sim 1$, where $k$
is the number of parallel copies.  Since the volume of $E'$ is
at least $A(n)^{-n} M^{-1}$, we also know that $k \lesssim A(n)^n M$.

Now suppose that $a^{k}$ does not vanish on $L(E)$.  Using the
multiple bisection lemma from Section 4, we can pick a
cycle $Z$ in $L(E)$ so that $Z$ bisects each set $U_i$.  Next we
investigate the directional volumes of a surface bisecting a copy
of $E'$.

Suppose that $Z$ bisects $E'$.  Let $E'_1, ..., E'_n$ be the
lengths of the principal axes of $E'$.  Let $e_1, ..., e_n$ be
unit length vectors with $e_j$ lying on the $j^{th}$ principal
axis of $E'$.  (To check the notation, each point $\pm E'_j e_j$
lies on the boundary of $E'$.)  The vectors $e_1, ..., e_n$ form an
orthonormal basis of $\mathbb{R}^n$.

\begin{lemma} Under the hypotheses in the last paragraph, the
following estimate holds for some $1 \le j \le n$:

$$\bar V_{Z \cap E'} (e_j) \gtrsim Vol(E') / E_j'.$$

\end{lemma}

\proof Let $L$ be a linear map taking $E'$ diffeomorphically to
the unit ball.  The map $L$ is diagonal with respect to the basis
$e_j$: in this basis, it scales the $j^{th}$ coordinate by
$1/E'_j$.  Then $L(Z)$ bisects the unit ball.  According to the
bisection clause in Lemma 4.1, $\bar V_{L(Z) \cap B} (e_j)
\gtrsim 1$ for some $j$.  When we change coordinates back and
interpret this inequality in $E'$, it gives the lemma.  We now
explain the coordinate change in detail.  We let $\pi_j$ denote
the orthogonal projection from $\mathbb{R}^n$ to $e_j^\perp$. 
Next we use equation 2 from section 2 to write directional
volumes in terms of $\pi_j$:

$V_{L(Z') \cap B}(e_j) = \int_{e_j^\perp} |L(Z') \cap B \cap
\pi_j^{-1}(y)| dy.$

$V_{Z' \cap E}(e_j) = \int_{e_j^\perp} |Z' \cap E \cap
\pi_j^{-1}(y)| dy.$

Comparing the right-hand sides we get the following formula:

$V_{Z' \cap E}(e_j) = (\prod_{i=1}^n E_i') (1 / E_j') V_{L(Z') \cap
B} (e_j).$

Averaging over $Z'$, we get an inequality for the mollified
directional volumes:

$\bar V_{Z \cap E}(e_j) = (\prod_{i=1}^n E_i') (1 / E_j') \bar
V_{L(Z) \cap B}(e_j) \gtrsim Vol(E') / E_j'.$

\endproof

Since $U_i$ is a translation of $E'$, we get the following
estimate:

For each $i$, there is some coordinate $j$, so that $\bar V_{Z
\cap U_i} (e_j) \gtrsim Vol(E') [E'_j]^{-1}$.

The number of translated ellipsoids $U_i$ is $k$, where $Vol(E')
k \sim 1$.  Combining our last estimate over all these
ellipsoids, we see that for a popular coordinate $j$, $\bar
V_{Z \cap B(p,1)} (e_j) \gtrsim (E'_j)^{-1} = A(n) E_j^{-1}$.

Now we choose $A(n)$ sufficiently large compared to our
dimensional constants, and we conclude that $\bar V_{Z \cap
B(p,1)} (e_j) > E_j^{-1}$, and so $\bar V_{Z \cap B(p,1)} (E_j
e_j) > 1$.  But the vector $E_j e_j$ is contained in $E$.  By the
definition of $L(E)$, we should have $\bar V_{Z \cap B(p,1)}(v)
\le 1$ for every $v \in E$.  This contradiction shows that our
assumption was wrong, and $a^{k}$ vanishes on $L(E)$.  But $k
\lesssim A(n)^n M$, and so $a^{C(n) M}$ vanishes on $L(E)$ for an
appropriate dimensional constant $C(n)$.

Reinspecting the argument we see that $a^{C_n M}$ vanishes on the
union of $\mathbb{RP}^n - Bi(U_i)$.  This latter set is open and
we have shown that it contains $L(E)$, and so $a^{C_n M}$
vanishes on a neighborhood of $L(E)$. \endproof

Next we explain how to upgrade this weak visibility lemma to get
the visibility lemma we originally stated.  For each sufficiently
large ellipsoid $E$, we have seen that $a^{C_n M}$ vanishes on
$L(E)$.  Remarkably, $a^{C_n M}$ vanishes on the union $\cup_E
L(E)$ as $E$ varies over all ellipsoids with volume at least
$M^{-1}$.  We can use the vanishing lemma to show that $a^{p
C_n M}$ vanishes on the union of any $p$ sets $L(E_k)$, but we
don't have any good control of the size of $p$.  The situation is
analogous to the following proposition, which is used in Gromov's
paper \cite{Gr2}.

\begin{refprop} (Gromov) Let $X$ be a manifold and let 
$f: X \rightarrow \mathbb{R}^m$ be a
map.  Suppose that for each unit ball $B(y,1)$ in $\mathbb{R}^m$,
the cohomology class $\alpha \in H^*(X)$ vanishes on
$f^{-1}[B(y,1)]$.  Then $\alpha^{m+1}$ vanishes on all of $X$.
\end{refprop}

\proof Triangulate $\mathbb{R}^m$ so that each simplex has
diameter at most 1/4.  Let $U_i$ be an open cover on
$\mathbb{R}^m$, indexed by the simplices of the triangulation
(including simplices of all dimensions).  It is possible to
choose $U_i$ in such a way that $U_i$ intersects $U_j$ only if
one of the corresponding simplices contains the other one.  In
particular, the open sets corresponding to two simplices of the
same dimension never intersect.  Also, each $U_i$ is contained in
a (1/10)-neighborhood of the corresponding simplex.  Since each
$U_i$ has diameter at most $1/2$, each $U_i$ is contained in some
unit ball, and so $\alpha$ vanishes on $f^{-1}(U_i)$.  Now for $0
\le l \le m$, let $V_l$ denote the union of $U_i$ as $i$ varies
among all the l-dimensional simplices of our triangulation.  For
any two l-dimensional simplices, $i_1$ and $i_2$, the
corresponding sets $U_{i_1}$ and $U_{i_2}$ are disjoint, and
hence their preimages $f^{-1}(U_{i_1})$ and $f^{-1}(U_{i_2})$ are
disjoint open subsets of $X$.  Therefore, $\alpha$ vanishes on
$f^{-1}(V_l)$.  Finally, by Lusternik-Schnirelmann,
$\alpha^{m+1}$ vanishes on all of $X$. \endproof

(The clever covering for $\mathbb{R}^m$ that appears here
originated in dimension theory.  See the book \cite{D} for
more information.)

Our argument is a variation on the proof of this proposition. 
The role of the space $\mathbb{R}^m$ is played by the space of
all ellipsoids in $\mathbb{R}^n$.  

Let $Ell$ denote the set of all closed ellipsoids in
$\mathbb{R}^n$ centered at the origin.  We put a distance
function on $Ell$ by saying that $dist_{Ell}(E_1, E_2) \le \log
D$ iff $(1/D) E_1 \subset E_2 \subset D E_1$.  We let
$Ell[M]$ denote the set of all ellipsoids contained in the unit
ball with volume at least $M^{-1}$.  Then we choose a maximal
1-separated subset of $Ell[M]$, given by finitely many ellipsoids
$Ell_1, ..., Ell_s$. The number of ellipsoids is finite, but it
grows exponentially with $M$.

Recall that $L_M$ is the set of hypersurfaces $Z \in \mathbb{RP}^N$ so
that $\bar Vis[ Z \cap B(p,1)] \le M$.  Next we divide the set $L_M$ into classes.  For any hypersurface
$Z$, we let $K[Z]$ be the convex set $\{ v \textrm{ such that }
|v| \le 1 \textrm{ and } \bar V_{Z \cap B(p,1)}(v) \le 1 \}$. 
We say that a hypersurface $Z$ lies in $A_k$ iff $K[Z]$ resembles
$Ell_k$ in the sense that $(10 n)^{-1/2} Ell_k \subset K[Z] \subset
(10n)^{1/2} Ell_k$.  Because our mollified function $\bar V_{Z \cap
B(p,1)}$ is continuous, the sets $A_k \subset \mathbb{RP}^N$ are
closed.

According to a lemma of Fritz John, any symmetric convex set $K$
can be approximated by an ellipsoid $E$ in the sense that
$n^{-1/2} E \subset K \subset n^{1/2} E$.  From this estimate, it follows
that the sets $A_k$ cover $L_M$.  

On the other hand, $A_k \subset L[(10 n)^{-1/2} Ell_k]$.  By the
weak visibility lemma, we see that $a^{C_n M}$ vanishes on a
neighborhood of each $A_k$.

Two sets $A_k$ and $A_l$ overlap only if the corresponding
ellipsoids $Ell_k$ and $Ell_l$ lie within a distance $C(n)$ of
each other, using our metric on $Ell$.  We want to bound
the multiplicity of the cover of $L_M$ by the sets $A_k$.  It suffices
to bound the number of ellipsoids $Ell_k$ inside a ball of radius $C(n)$
in the space $Ell$.  Let $Ell_0$ denote the unit ball.  The closed ball of
radius $C(n)$ around $Ell_0$ is a compact subset of $Ell$.  (The space
$Ell$ is a finite-dimensional manifold, and our metric defines the usual
topology on the manifold.)  By compactness, any set of 1-separated ellipsoids
$Ell_i$ inside this ball has cardinality bounded by some $C'(n)$.  But
there is nothing special about the unit ball $Ell_0$.  In fact, the space $Ell$
is extremely symmetrical.  The group $GL(n, \mathbb{R})$ acts on $Ell$
in the following way.  Given a linear map $M \in GL(n, \mathbb{R})$ and
an ellipsoid $E \in Ell$, we define $M(E)$ to be the image of $E$
under the map $M$.  This group action is an isometry using our metric
on $Ell$.  It is also transitive because of the spectral theorem.  Therefore,
any ball of radius $C(n)$ contains at most $C'(n)$ 1-separated points.  Hence
the multiplicity of the cover $\{ A_k \}$ is bounded by $C'(n)$.

Let $B_k$ be tiny open neighborhoods of $A_k$ so that $a^{C_n M}$
vanishes on $B_k$.  Since the sets $A_k$ are closed, we can
arrange that $B_k$ and $B_l$ intersect only if $A_k$ and $A_l$
intersect.  The $B_k$ form an open cover of a neighborhood of
$L_M$ with multiplicity at most $C'(n)$.  We color the sets $B_k$
using $C'(n)$ colors so that overlapping sets have distinct
colors.  For each color $\alpha$ from $1$ to $C'(n)$, we let
$C_\alpha$ denote the union of all sets $B_k$ with the color
$\alpha$.  Because these sets are disjoint, $a^{C_n M}$ vanishes
on $C_\alpha$ for each $\alpha$.  Now by the
Lusternik-Schnirelmann vanishing lemma, $a^{C''_n M}$ vanishes on
the union of $C_\alpha$, which includes a neighborhood of $L_M$.
\endproof

Combining the visibility lemma and the Lusternik-Schnirelmann
vanishing lemma, we can find a degree $d$ algebraic hypersurface
with large visibility on various cubes.  The following lemma is the
main result of this section.

\newtheorem*{lvmc}{Large visibility on many cubes}

\begin{lvmc} Consider the
standard unit lattice in $\mathbb{R}^n$.  Let $M$ be a function
from the set of n-cubes in the unit lattice to the non-negative
integers.  Then we can find an algebraic hypersurface of degree
$d$ so that $\overline{Vis}[Z \cap Q_k] \ge M(Q_k)$ for every cube
$Q_k$, where the degree $d$ is bounded by $C(n) [\sum_k
M(Q_k)]^{1/n}$.
\end{lvmc}

\proof The space of degree $d$ hypersurfaces is parametrized by
$\mathbb{RP}^N$, where $N = {{n + d} \choose d} - 1 \ge c(n)
d^n$. We let $a$ denote the fundamental cohomology class of
$\mathbb{RP}^N$. Let $S[Q_k]$ denote the set of surfaces $Z$
where $\overline{Vis}[Z \cap Q_k] < M(Q_k)$.  According to the
visibility lemma, the cohomology class $a^{C(n) M(Q_k)}$ vanishes
on a neighborhood of $S[Q_k]$.  By the Lusternik-Schnirelmann
vanishing lemma, the cohomology class $a^{C(n) \sum M(Q_k)}$
vanishes on a neighborhood of $\cup_k S[Q_k]$. But $a^N$ does not
vanish on $\mathbb{RP}^N$.  As long as $C(n) \sum M(Q_k) < c(n)
d^n \le N$, there is a variety $Z$ which does not lie in any
$S(Q_k)$.  Unwinding the definition, we see that $\overline{Vis}[Z
\cap Q_k] \ge M(Q_k)$ for every $k$.  Our condition on $d$ is
$C(n) \sum M(Q_k) < c(n) d^n$, which holds for any $d > C'(n)
[\sum_k M(Q_k)]^{1/n}$.  \endproof

\section{Multilinear Kakeya estimates}

Let us recall the setting of the multilinear Kakeya estimate.  We have some
unit cylinders $T_{j,a}$, where $1 \le j \le n$, while
for each $j$, $1 \le a \le A(j)$.  We let $v_{j,a}$ be a unit
vector parallel to the core of $T_{j,a}$.  We assume that
tubes with different values of $j$ are quantitatively transverse in the following
sense.  For any sequence of tubes $T_{1,
a(1)}$, ..., $T_{n, a(n)}$, $|Det(v_{1, a(1)}, ... v_{n, a(n)})| \ge \theta.$

\begin{theorem} (Multilinear Kakeya estimate) Under the hypotheses in the last
paragraph, the following inequality holds.

$$\int \left[ \prod_{j=1}^n \left( \sum_{a=1}^{A(j)}
\chi_{T_{j,a}} \right) \right]^{\frac{1}{n-1}} < C(n) \theta^{-\frac{1}{n-1}}
\prod_{j=1}^n A(j)^{\frac{1}{n-1}}$$.  

\end{theorem}

\proof We consider the standard unit cube lattice in
$\mathbb{R}^n$.  For each cube $Q_k$ in this lattice, we define
the following functions, measuring how many tubes of different
types go through $Q_k$.  We let $M_j(Q_k)$ denote the number of tubes 
$T_{j,a}$ which go through $Q_k$.  Then we let $F(Q_k)$ be the product of these:

$$ F(Q_k) := \prod_{j=1}^n M_j(Q_k).$$

It suffices to prove the following estimate for $F(Q_k)$:

$$\sum_k F(Q_k)^{\frac{1}{n-1}} < C(n) \theta^{-\frac{1}{n-1}}
\prod_{j=1}^n A(j)^{\frac{1}{n-1}}. \eqno{(*)}$$

Since we have only finitely many tubes, the function $F(Q_k)$
vanishes outside of finitely many cubes.  We fix a large cube of
side length $S$ containing all of the relevant cubes $Q_k$.  Next we apply
the Large-visibility-on-many-cubes lemma from Section 5.
The lemma guarantees that we can find a hypersurface $Z_0$ of
degree $d \lesssim S$ obeying the following visibility estimates.  For
every cube $Q_k$, 

$$\overline{Vis}[Z_0 \cap Q_k] \ge S^n F(Q_k)^{\frac{1}{n-1}} \left[ \sum_k
F(Q_k)^{\frac{1}{n-1}} \right]^{-1}. \eqno{(1)}$$

Adding $C(n) S$ hyperplanes to $Z_0$, we produce a variety $Z$ of
degree $d \lesssim S$ that still obeys the visibility estimate
above, and also $\bar V_{Z \cap Q_k}(v) \ge |v|$ for each cube
$Q_k$ where $F(Q_k) > 0$.  Equation 1 gives a strong lower bound for $V_{Z \cap Q_k}(v)$
in some directions, and this last estimate gives a weak lower bound in all
directions.

Next we apply the cylinder estimate to control the directed
volumes of $Z$ in cubes along a given tube $T_{j,a}$.  (The
estimate we need is the cylinder clause of Lemma 5.1.)  For each
tube $T_{j,a}$, we have the following estimate:

$$\sum_{Q_k \textrm{ that intersect } T_{j,a}} \bar V_{Z \cap
Q_k} (v_{j,a}) \lesssim S. $$

We would like to sum this inequality over all $a$ from $1$ to
$A(j)$, but the vectors $v_{j,a}$ are changing.

For each cube $Q_k$ and each $j$, we pick a vector $v_{j,k}$ from
among $v_{j,a}$ so that $\bar V_{Z \cap Q_k} (v_{j,k}) =
\min_{a=1}^{A(j)} \bar V_{Z \cap Q_k}(v_{j,a})$.

Substituting $v_{j,k}$ for $v_{j,a}$ in the last inequality and
summing over $a$ yields the following.

$$\sum_{Q_k} M_j(Q_k) \bar V_{Z \cap Q_k} (v_{j,k}) \lesssim S
A(j). \eqno{(2)}$$

Next we need a lemma relating $\overline{Vis}$ and $\bar V$.

\begin{lemma} For each cube $Q_k$, the following inequality
holds.

$$ \overline{Vis}[Z \cap Q_k] < C(n) \theta^{-1} \prod_{j=1}^n \bar
V_{Z \cap Q_k} (v_{j,k}). $$

\end{lemma}

\proof Let $v'_{j,k} = v_{j,k} / [ \bar V_{Z \cap Q_k}(v_{j,k})]$. 
Because we added the hyperplanes to $Z_0$, we know that
$\bar V_{Z \cap Q_k}(v) \ge 1$ for all unit vectors $v$.  Hence
$|v'_{j,k}| \le 1$.  We know that $\bar V_{Z \cap
Q_k}(\pm v'_{j,k}) = 1$ for each $j$.  Since the directed volume is a
convex function of $v$, $\bar V_{Z \cap Q_k}(v) \le 1$ for
every $v$ in the convex hull of the $2n$ points $\pm v'_{j,k}$. 
This convex hull is contained in the unit ball, and its volume is
$c(n) det(v'_{1,k}, ..., v'_{n,k}) = c(n) [\prod_{j=1}^n \bar
V_{Z \cap Q_k}(v_{j,k})]^{-1} det (v_{1,k}, ..., v_{n,k})$. 
Because of our transverality assumption, the volume is $\ge
c(n) \theta [\prod_{j=1}^n \bar V_{Z \cap Q_k}(v_{j,k})]^{-1}$.
Hence the set of all $v$ with $|v| \le 1$ and $\bar V_{Z \cap Q_k}(v)
\le 1$ has volume at least $\gtrsim [\prod_{j=1}^n \bar
V_{Z \cap Q_k}(v_{j,k})]^{-1} \theta$.  The visibility $\bar
Vis[Z \cap Q_k]$ is the inverse of this volume, which is at most
$C(n) \theta^{-1} \prod_{j=1}^n \bar V_{Z \cap Q_k}(v_{j,k})$.
\endproof

Now we follow a string of inequalities powered by the visibility
estimate in Equation 1 and the cylinder estimate in Equation 2.

$$S \left[ \sum_k F(Q_k)^{\frac{1}{n-1}} \right] ^{\frac{n-1}{n}}
= \sum_k \left( S F(Q_k)^{\frac{1}{n-1}} / \left[ \sum_k
F(Q_k)^{\frac{1}{n-1}} \right]^{\frac{1}{n}} \right),$$

\noindent We use the visibility estimate in equation 1.

$$\le \sum_k F(Q_k)^{1/n} \overline{Vis}[Z \cap Q_k]^{1/n},$$

\noindent Now we apply H\"older to the products $F(Q_k) =
\prod_{j=1}^n M_j(Q_k)$ and $\overline{Vis}[Z \cap Q_k] \lesssim
\theta^{-1} \prod_{j=1}^n \bar V_{Z \cap Q_k} (v_{j,k})$.

$$\lesssim \theta^{-1/n} \prod_{j=1}^n \left[ \sum_k M_j(Q_k)
\bar V_{Z \cap Q_k} (v_{j,k}) \right]^{1/n},$$

\noindent We use the cylinder estimate in equation 2.

$$\lesssim \theta^{-1/n} \prod_{j=1}^n [S A(j)]^{1/n} = S
\theta^{-1/n} \prod_{j=1}^n A(j)^{1/n}. $$

\noindent We summarize this string of inequalities.

$$S \left[ \sum_k F(Q_k)^{\frac{1}{n-1}} \right]^{\frac{n-1}{n}} \lesssim S
\theta^{-1/n} \prod_{j=1}^n A(j)^{1/n}. $$

\noindent Finally, we cancel the $S$ on each side and raise the equation to
the power $\frac{n}{n-1}$.

$$\sum_k F(Q_k)^\frac{1}{n-1} \lesssim
\theta^{-\frac{1}{n-1}} \prod_{j=1}^n A(j)^{\frac{1}{n-1}}.$$

This establishes the inequality $(*)$ and hence the theorem.

\endproof

\section{Box estimates for unions of tubes}

The multilinear Kakeya estimate of Bennett, Carbery, and Tao
implies that Kakeya sets must be rather ``plany''.  Here we give
a quantitative estimate of planiness.

\begin{boxest} There is a constant $C(n) > 0$ so that the
following holds.  Suppose $X \subset \mathbb{R}^n$ is a union of
cylinders with radius 1 and length $L >> 1$.  For each $x \in X$
we can choose a convex set $B(x)$ with the following
properties.

1. The set $B(x)$ contains $x$.  In fact, $B(x)$ is a symmetric
convex body translated so that the center is $x$.  

2. The set $B(x)$ has volume at most $Vol(X)$.

3. For every cylinder $T \subset X$ of radius 1 and length $L$,
if we pick a random point $x \in T$, then the tube $T$ lies in
the rescaled set $\sigma B(x)$ with probability at least $1 -
C(n) \sigma^{-1}$.  (This probability estimate holds for every
$\sigma > C(n)^{-1}$.)

\end{boxest}

\proof We pick a collection of disjoint balls $B_i$ of radius
(1/10) so that the union of $3 B_i$ covers $X$.  The number of
balls is $\lesssim Vol(X)$.

We can assume that $Vol(X)$ is significantly less than $L^n$,
because otherwise we just take each $B(x)$ to be a cube with side
length $L$.  By the large-visibility-on-many-cubes lemma from
Section 5, we can choose an algebraic
hypersurface $Z$ so that $\overline{Vis}[Z \cap B_i] \ge L^n /
Vol(X)$ for each ball in our cover, with degree $\lesssim L$.

We use the hypersurface $Z$ to define our box function $B(x)$. 
First take the set $\textrm{ } \{ v \textrm{ such that } |v| \le
1 \textrm { and } \bar V_{Z \cap B(x,1)}(v) \le 1 \}$.  Let
$B_0(x)$ be the translate of this set so that it is centered at
$x$ instead of at the origin.  Then let $B(x)$ be the rescaling
of $B_0(x)$ by a factor $L$, keeping it centered at $x$.  For
each $x \in X$, the
unit ball $B(x,1)$ contains at least one ball $B_i$ from our set
of balls, and so $\overline{Vis}[Z \cap B(x,1)] \ge L^n / Vol(X)$. 
Therefore, the convex set $B_0(x)$ has volume at most $Vol(X) /
L^n$, and so the box $B(x)$ has volume at most $Vol(X)$.

Now fix a number $\sigma > 1$ and a tube $T \subset X$ with
radius 1 and length $L$.  Let $v$ be a unit vector pointing
parallel to the core of $T$.  If $x
\in T$, then $T$ lies in $\sigma B(x)$ unless $\bar V_{Z \cap
B(x,1)}(v) \ge \sigma/2$.

On the other hand, we will estimate the average value of $\bar
V_{Z \cap B(x, 1)}(v)$ as $x$ varies in $T$.  

\begin{lemma} Let $Z'$ denote any algebraic
hypersurface of degree $\lesssim L$.  Then the average value of
$V_{Z' \cap B(x,1)}(v)$ over $x$ in $T$ is bounded as follows.

$$|T|^{-1} \int_{T} V_{Z' \cap B(x,1)}(v) dx \lesssim 1.$$

\end{lemma}

This lemma is essentially the cylinder estimate Lemma 2.1, as we
will see below.  Given the lemma, we can finish the proof of the
box estimate.  Applying the lemma to averages over appropriate
$Z'$, we get the following estimate for the mollified volume
$\bar V$:

$$|T|^{-1} \int_T \bar V_{Z \cap B(x,1)}(v) dx \lesssim 1.$$

Let $B \subset T$ be the set of bad points where $T$ is not
contained in $\sigma B(x)$.  At each bad point, $\bar V_{Z
\cap B(x,1)}(v) \ge \sigma/2$.  Since the average value of $\bar
V_{Z \cap B(x,1)}(v)$ is at most $C(n)$, it follows that the
volume of $B$ is at most $2 C(n) \sigma^{-1} |T|$. \endproof

Now we turn to the proof of Lemma 7.1.

\proof We want to understand the average: $|T|^{-1} \int_T V_{Z'
\cap B(x,1)}(v) dx$.  The directional volume is itself an integral. 
We expand that integral and apply Fubini:

$$|T|^{-1} \int_T V_{Z' \cap
B(x,1)}(v) dx = |T|^{-1} \int_T \left( \int_{Z' \cap B(x,1)} |N(y)
\cdot v| dy \right) dx \le $$

$$|T|^{-1} \int_{Z' \cap 3T} |N(y) \cdot v| \left( \int_{B(y,1)}
dx \right) dy
= C(n) L^{-1} V_{Z' \cap 3T}(v).$$

But according to the cylinder estimate Lemma 2.1, $V_{Z' \cap
3T}(v) \lesssim L$.  Plugging in, we see that our average is
$\lesssim 1$. \endproof

\section{Appendix: The Lusternik-Schnirelmann vanishing lemma}

In this section, we give a proof of the vanishing lemma.  There are proofs in the literature, but I
will try to write the proof here in a way that's accessible with a minimum of background in algebraic topology.

First, I will prove the lemma in the special case of de Rham cohomology on a manifold.  This setting is probably familiar to more readers, and the proof in this setting is clearest.  In the paper, we have to apply the vanishing lemma to mod 2 cohomology, so we do the general case afterwards.

\newtheorem*{vanlemma1}{Vanishing lemma for de Rham cohomology}

\begin{vanlemma1} (not optimal version) Let $M$ be a smooth manifold.  Let $a_1$ and $a_2$ be cohomology classes
in $H^*(M, \mathbb{R})$.  Suppose that $a_1 $ vanishes on some open set $S_1 \subset X$ and that $a_2$ vanishes on some open set $S_2 \subset X$.  Let $K \subset S_1 \cup S_2$ be a compact set Then the cup product $a_1 \cup
a_2$ vanishes on $K$.
\end{vanlemma1}

In fact, $a_1 \cup a_2$ vanishes on all of $S_1 \cup S_2$, not just on the compact subsets.  But I chose to prove the weaker statement above because it makes the proof shorter and clearer.

\begin{proof}  Because we are using cohomology with real coefficients and working on a manifold, we may use de Rham cohomology.  Let $\alpha_1$ be a differential form that represents the cohomology class $a_1$.  The first point of the proof is that we can choose $\alpha_1$ to vanish on almost all of $S_1$.  Let's see how to do this.  We know that the restriction of $a_1$ to $S_1 \subset M$ is zero.  In other words, the
restriction of $\alpha_1$ to $S_1$ is exact.  In other words, there is a form $\beta$ on $S_1$ so that $d \beta = \alpha_1$ on $S_1$.  The form $\beta$ is only defined on $S_1$.  Now let $K_1 \subset S_1$ be any compact
subset - the reader should imagine that $K_1$ is almost all of $S_1$.  We can find a form $\beta'$ on all of $M$ so that
$\beta'$ restricted to $K_1$ agrees with $\beta$.  Hence $d \beta'$ is an exact form on 
all of $M$.  Also $d \beta' = \alpha_1 $ on $K_1$.  Since $d \beta'$ is exact, $\alpha_1 - d \beta'$ still represents the cohomology class $a_1$.  But $\alpha_1 - d \beta'$ vanishes pointwise on $K_1$.

By the previous paragraph, we may pick a differential form $\alpha_1$ on $M$ which represents $a_1$ and vanishes pointwise on $K_1$.  By the same argument, for any compact $K_2 \subset S_2$, we can
pick a differential form $\alpha_2$ on $M$ which represents $a_2$ and vanishes pointwise on $K_2$.  
Now the wedge product of forms $\alpha_1 \wedge \alpha_2$ represents the cup product $a_1 \cup a_2$.  On the other hand, the wedge product $\alpha_1 \wedge \alpha_2$ vanishes pointwise on $K_1 \cup K_2$.  Hence $a_1 \cup a_2$ vanishes on $K_1 \cup K_2$.  Therefore, $a_1 \cup a_2$ vanishes on any compact subset $K \subset S_1 \cup S_2$.  \end{proof}

We now prove the vanishing lemma in general.

\begin{vanlemma} Let $X$ be a CW complex - for example, a manifold.  Let $R$ be any ring.  Let $a_1, a_2$ be cohomology classes in $H^*(X, R)$.  Suppose that $a_1 $ vanishes on some open set $S_1 \subset X$ and that $a_2$ vanishes on some open set $S_2 \subset X$.  Then the cup product $a_1 \cup a_2 $ vanishes on the union $S_1 \cup S_2$.
\end{vanlemma}

\begin{proof}  This time we work with singular cohomology.  Singular cohomology and cup products are well explained in Hatcher's book on algebraic topology, \cite{H}, chapters 3.1 and 3.2.   Let $\alpha_1$ be a singular cocycle representing $a_1$.  We know that $a_1$ restricted to $S_1$ is zero.  Therefore, we can choose a singular cochain $\beta$ on $S_1$ so that $\partial \beta$ is equal to the restriction of $\alpha_1$ to $S_1$.  We can automatically extend $\beta$ to a singular cochain $\beta'$  on all of $X$.
Then we look at the cocycle $\alpha_1 - \partial \beta'$.  Since $\partial \beta'$ is exact, this cocycle
still represents the cohomology class $a_1$.  The cocycle $\alpha_1 - \partial \beta'$ vanishes on
any chain supported in $S_1$.

By the previous paragraph, we may pick a singular cocycle $\alpha_1$ representing $a_1$ so that $\alpha_1$ vanishes on $S_1$.  Similarly, we may pick a cocycle $\alpha_2$ representing $a_2$ so
that $\alpha_2$ vanishes on $S_2$.  Now we look at the cup product $\alpha_1 \cup \alpha_2$,
which represents $a_1 \cup a_2$.  The product $\alpha_1 \cup \alpha_2$ vanishes on any singular
simplex supported in $S_1$ or supported in $S_2$.

Now let $f$ denote a singular simplex supported in $S_1 \cup S_2$.  We will subdivide $S$ into small pieces so that each piece lies in either $S_1$ or $S_2$.  Recall that $f$ is a continuous map from the simplex $\Delta$ to $X$.  We subdivide the simplex into many small simplices.  Restricting $f$ to each small simplex, we get various maps $g_i$ from $\Delta$ to $X$.  The sum $\sum g_i$ is a singular chain that parametrizes the image of $f$.  Now it's not true that $f = \sum g_i$ as singular chains.  But it is true
that $f - \sum g_i$ is a boundary.  Since $\alpha_1 \cup \alpha_2$ is a cocycle, $\alpha_1 \cup \alpha_2
(f) = \sum_i \alpha_1 \cup \alpha_2 (g_i)$.  If we subdivide finely enough, then each $g_i$ is contained
in either $S_1$ or $S_2$. (At this step, we use the fact that $S_1$ and $S_2$ are open.)   So each term $\alpha_1 \cup \alpha_2 (g_i)$ vanishes.  Hence $\alpha_1 \cup \alpha_2 (f) = 0$ for any singular simplex $f$ in $S_1 \cup S_2$.  In other words, $\alpha_1 \cup \alpha_2$ vanishes
on $S_1 \cup S_2$.  Hence the cup product $a_1 \cup a_2$ vanishes on $S_1 \cup S_2$. \end{proof}

\end{document}